\documentclass[11pt]{article}

\usepackage{amsfonts}
\usepackage{amssymb}
\usepackage{amsmath}
\usepackage{amsthm}
\usepackage{latexsym}
\usepackage{verbatim}
\usepackage{epsfig}
\usepackage{amsbsy}
\usepackage{amscd}
\usepackage{bm}

\usepackage{graphicx}

\usepackage{ifthen} 
\usepackage{xcolor}
\usepackage[colorlinks=true, linkcolor=teal, anchorcolor=teal,
            citecolor=blue,  filecolor=blue, menucolor=blue, pagecolor=teal,											urlcolor=blue, plainpages=false, pdfpagelabels]{hyperref}

\newcommand{\mymail}[1]{\href{mailto:#1}{\texttt{#1}}}

\usepackage{fancyhdr}
\usepackage[realmainfile]{currfile}
\fancypagestyle{firststyle}{

\cfoot{{\footnotesize $^\ast$: Supported by XXX $\qquad \qquad$ $^\dagger$: Partially supported by NSF grants \#2120318 and \#2411270}} 
}

\newcommand{\setauthA}[1]{\def\authA{#1}}
\newcommand{\setauthB}[1]{\def\authB{#1}}
\newcommand{\setauthC}[1]{\def\authC{#1}}

\def\printA{\begin{tabular}{l} \authA \end{tabular}}
\def\printB{\begin{tabular}{l} \authB \end{tabular}}
\def\printC{\begin{tabular}{l} \authC \end{tabular}}

\newcommand{\makemytitle}[1]{\begin{center}{\textsf{\LARGE #1}}
  \end{center}
}

\usepackage[sf,bf]{titlesec}

\usepackage{algorithmic}
\usepackage{algorithm}

\usepackage[nonamebreak,square,numbers,sort]{natbib}

\usepackage[letterpaper, textheight = 676pt]{geometry}
\providecommand{\M}[1]{\mathbf#1}
\providecommand{\mc}[1]{\mathcal#1}
\providecommand{\mc}[1]{\mathcal#1}

\newcommand{\R}{{\mathbb R}}

\DeclareMathOperator{\E}{\mathbb{E}}

\providecommand{\T}{\top} 

\providecommand{\wt}[1]{\widetilde{#1}}

\providecommand{\nnorm}[1]{ \lVert#1 \rVert}

\newcommand{\nscp}[2]{\langle#1, #2\rangle}

\newcommand{\blanco}[1]{  }

\newcommand{\deriv}[3]{%
\ifthenelse{#1 = 1}{\frac{d\,#2}{d\,#3}}{\frac{d^{{#1}} #2}{d{#3}^{{#1}}}}
}

\newcommand{\partials}[3]{%
\ifthenelse{#1 = 1}{\frac{\partial\,#2}{\partial\,#3}}{\frac{\partial^{#1}
    #2}{\partial#3^{#1}}}
}

\def \coloneq{\mathrel{\mathop:}=}
\def \invcoloneq{=\mathrel{\mathop:}}
\def \eps{\varepsilon}

\newtheorem{theo}{Theorem}

\newtheorem{lemmachen}{Theorem}
\newtheorem{corollary}{Theorem}

\newtheorem{theorem}[theo]{Theorem}

\newtheorem{corro}[corollary]{Corollary}
 \newtheorem{lemma}[lemmachen]{Lemma}

\newenvironment{bew}{\begin{proof}[Proof]}{\end{proof}}

\def\R{\mathbb{R}}

\def\eps{\epsilon}

\usepackage{subfigure}
\usepackage{caption}
\usepackage{dsfont}
\usepackage{mathrsfs}

\newboolean{isjournal}
\setboolean{isjournal}{true}

\newboolean{nocolors}
\setboolean{nocolors}{false}

\newboolean{usemtpro}
\setboolean{usemtpro}{false}

\ifthenelse{\boolean{nocolors}}{\hypersetup{colorlinks=false}}{}

\makeatletter
\newcommand\footnoteref[1]{\protected@xdef\@thefnmark{\ref{#1}}\@footnotemark}
\makeatother

\setauthC{{\bfseries Jonathan Steffani$^{1}$}}

\setauthB{{\bfseries Martin Slawski$^{2\dagger}$}}

\setauthA{{\bfseries Fadoua Balabdaoui}$^{1*}$}

\begin{document}
\thispagestyle{firststyle}

\makemytitle{{\Large {\bfseries {Identifiability in Unlinked Linear Regression:\\[.5ex] Some Results and Open Problems}}}}
\vskip 3.5ex
%
{\large\begin{center}
\printA
\printB
\printC
\vskip1.5ex
{\scriptsize $^{1}$Department of Mathematics, ETH Z\"urich, 8092 Zurich, Switzerland $\; \; \;$}\\
{\scriptsize $^{2}$Department of Statistics, University of Virginia, Charlottesville, VA 22903, USA $\; \; \;$}\\[.5ex] 
{\footnotesize \mymail{fadoua.balabdaoui@stat.math.ethz.ch} 
$\quad$ \mymail{ebh3ep@virginia.edu} $\quad$ \mymail{jonathansteffani@web.de}} 
\end{center}}
\vskip 3.5ex


\begin{abstract} 
A tacit assumption in classical linear regression problems is the full knowledge of the existing link between the covariates and responses. In \textit{Unlinked linear regression} (ULR) this link is either partially or completely  missing. While the reasons causing such missingness can be different, a common challenge in statistical inference is the potential non-identifiability of the regression parameter. In this note, we review the existing literature on identifiability when the $d \ge 2$ components of the vector of covariates are independent and identically distributed. When these components have different distributions, we show that it is not possible to prove similar theorems in the general case. Nevertheless, we prove some identifiability results, either under additional parametric assumptions for $d \ge 2$ or conditions on the fourth moments in the case $d=2$. Finally, we draw some interesting connections between the ULR  and the well established field of Independent Component Analysis (ICA). 
\end{abstract}


\section{Introduction}\label{sec:intro}
\begin{figure}
\begin{center}
\begin{tabular}{cc}
$\begin{array}{c} 
\includegraphics[width = 0.3\textwidth]{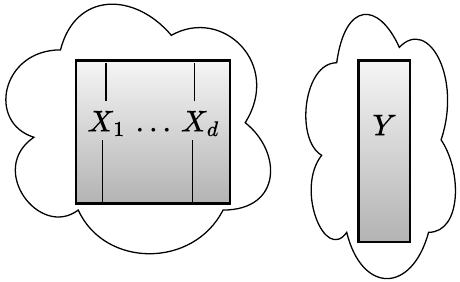} \\[23ex]
\end{array}$ \hspace*{2ex} & 
\hspace*{2ex} \includegraphics[width = 0.3\textwidth]{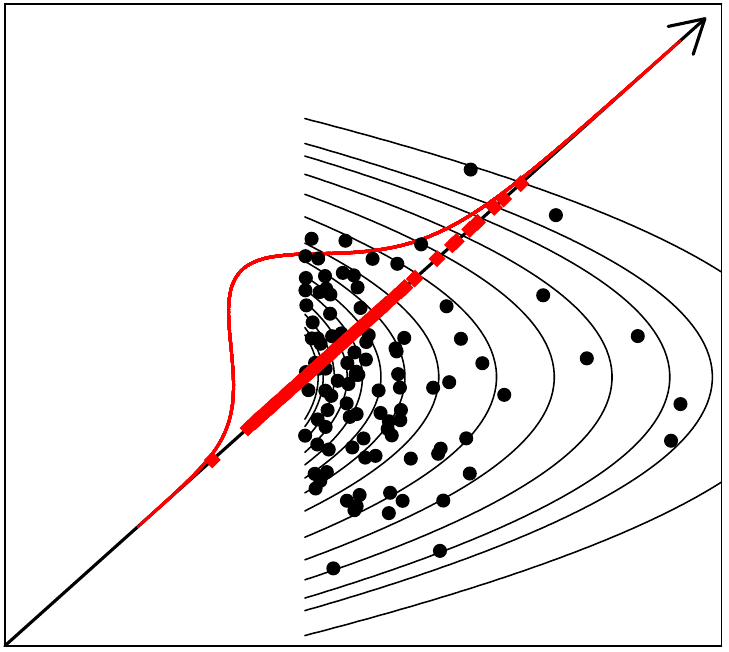}
\end{tabular}
\end{center}
\vspace*{-22ex}
\caption{Illustration of Unlinked Linear Regression. Left: The setup motivating the problem -- two data sources containing predictors and response, respectively, not necessarily of the same size nor necessarily corresponding to the same statistical units. Right: $X = (X_1, X_2)$ with $X_1 \sim \textrm{Exp}(1)$ and $X_2 \sim \mc{N}(0,1)$ and $Y \overset{d}{=} X_1 + X_2$. The contour lines of the density of $X$ are represented by solid black lines, and the black dots show 100 i.i.d.~realizations of $X$. The corresponding $Y$'s and their density (shown in color) are obtained via projection on the vector $(1,1)^{\T}$. As explained below, the situation depicted here is that of ``weak identifiability". }\label{fig:unlinked_illu}
\end{figure}

Linear regression is a basic tool in every data analyst's toolbox. A tacit assumption is that the predictor
variables $X = (X_1,\ldots,X_d)$ and the response variable $Y$ are observed concurrently for each statistical unit under consideration. Departures from this assumption are studied in (i) the missing data literature with an emphasis of dealing with various forms of missingness of (subsets of the) predictors or responses for some of the observations \cite[e.g.,][]{LittleRubin2020, KimShao2021}, and (ii) in semi-supervised learning scenarios involving missing responses \cite{Chapelle2006}. A third setting that has not attracted notable interest until recently concerns situations in which predictors and responses reside
in different data sources as depicted in Figure \ref{fig:unlinked_illu}. 

A common solution is to unify the two files by linking them, a process that is typically referred to as ({\em Record}) {\em Linkage} (RL, \citep[e.g.][]{Newcombe, Binette2022, Christen2012, Winkler2010}). While simple to conceptualize, it is well-known that in practice, RL tends to be problem-ridden. Often there are no unique identifiers available that would enable the straightforward identification of matching records; the use of quasi-identifiers such as names, dates, addresses, etc.~render the process laborious and ad-hoc, frequently relying on heuristics to determine matches and non-matches \cite{abowd2021}. Dealing with the resulting uncertainty in downstream analysis such as regression can be a challenging task \cite{Slawski2024, Kamat2024}. 

In addition, RL can be computationally intensive even with state-of-the-art implementations \cite{Enamorado2018}, raises concerns from the point of view privacy \cite{Sweeney2001, Narayanan2008}, and may not be reconcilable with reproducibility standards: if quasi-identifiers contain personally identifiable information (PII) details of the linkage cannot be published. Even when ignoring the privacy leakage induced by quasi-identifiers, it is well-known that the piling of information at a single location via file linkage may greatly increase re-identification risk (i.e., being able to associate an entry in database
with a specific individual). Indeed, disaggregating a database containing many variables into multiple smaller slices is an established method of statistical disclosure control \cite{Paass1988}. 

While the above points already provide ample motivation for studying {\em Unlinked Linear Regression} (or more generally, {\em Unlinked Data Analysis}), i.e., roughly speaking, the joint analysis of multiple sets of variables without an underlying integrated set of microdata, Unlinked Data Analysis extends beyond the realm 
of analyses that can be performed using RL. The latter requires that the statistical units underlying the two files are identical or at least exhibit substantial overlap. However, this may be a serious restriction: consider, e.g., the situation in which the two data sources originate from two surveys pertaining to the same population, in which case it may still be possible to perform inference regarding the association of variables across files. For instance, the approach adopted in Ecological Regression \cite{King1997, Gelman2001} is to compute averages within common strata and pretend that these averages constitute pairwise microdata. 

Unlinked Regression as first formulated in \cite{carp16, balabdaoui2021unlinked} for $d = 1$ hinges on the idea that the underlying random variable $Y$ is a function of the random variable $X$, possibly subject to additive noise. Several other papers \cite{rigollet2019uncoupled, Meis2021, durot2024minimax} have studied the case in which the function is monotone (with known direction of monotonicity), including a generalization to $d \geq 1$ based on the notion of transport maps \cite{SlawskiSen2022}. The methodology developed in these papers is tightly linked to deconvolution problems \cite{meister2009} and techniques employed in that context. As a result, unlinked monotone regression tends to exhibit similar rates of 
convergence, which can be dramatically slower than in the conventional setup with paired data.

Unlinked Linear Regression (ULR) has been studied recently in \cite{monafadoua}. The problem stated therein generalizes {\em Permuted} or {\em Shuffled} Linear Regression in which the two files are of the same size and based on the same statistical units. The latter problem has garnered much interest in signal processing \cite{Unnikrishnan:2018gp, Pananjady:2018hd}, machine learning \cite{Hsu2017, TsakirisICML19}, and statistics \cite{SlawskiBenDavid2017}. In brief, ULR postulates that 
\begin{eqnarray}\label{ULR}
Y \overset{d}{=} \beta_0^{\T} X + \eps,
\end{eqnarray}
where $\overset{d}{=}$ denotes equality in distribution, $\beta_0 \in \R^d$ is a vector of regression coefficients and $\eps$ is additive noise independent of the random vector $X$. A visual representation illustrating this setup for a specific example (without noise) is provided in Figure \ref{fig:unlinked_illu}. In geometrical terms, ULR aims to recover the direction of the projection from the distribution of the projection and the distribution of the input variables. 

As in classical linear regression, the primary goal of ULR is to make inference about the unknown regression parameter. Clearly, this task can only be accomplished if identifiability holds, i.e., if the regression 
parameter is uniquely determined in a suitable sense. Plain identifiability means that 
\begin{eqnarray}\label{eq:plain_identifiability}
\beta_0^{\T}  X  \stackrel{d}{=} \beta_1^{\T
}X    \Longleftrightarrow  \beta_0  = \beta_1, 
\end{eqnarray}
with weakened notions to be considered below. 

In the monotone regression case mentioned above, 
identifiability is established in \cite{balabdaoui2021unlinked}. In their work on ULR, 
\cite{monafadoua} establish guarantees such as asymptotic Normality for their Deconvolution Least Squares estimator (DLSE) {\em assuming identifiability}. A specific positive example provided therein is as follows:  
let $X= (X_1, X_2)^{\T}$ such that $X_1 \sim \mathcal{N}(\mu, \sigma^2)$ is independent $X_2 \sim \mathrm{Exp}(\lambda)$ for $\mu \ne 0$, $\sigma > 0$ and  $\lambda > 0$, then the ULR model
\eqref{ULR} is identifiable.  

On the other hand, it is shown that the same model is non-identifiable if $X \sim \mathcal{N}(\mu, \Sigma)$ for any $\mu \in \mathbb R^d$ and $\Sigma$.  To see this, let us assume that $\mu = \mathbf{0} \in \mathbb R^d$ and $\Sigma$ is positive definite.  If $\beta_0$ satisfies the model \eqref{ULR} then the same holds for any $\beta \in \R^d$ such that
\begin{eqnarray*}
\Vert \Sigma^{1/2} \beta  \Vert_2 = \Vert \Sigma^{1/2} \beta_0  \Vert_2 
\end{eqnarray*}
where $\nnorm{\cdot}_2$ denotes the Euclidean norm.  In particular, this implies that when $X$ is a standard Gaussian, any vector in the sphere of radius $\nnorm{\beta_0}_2$ is an adequate regression vector.  In their Theorem 7, the authors of \cite{monafadoua} show that if $X = (X_1, \ldots, X_d)^{\T}$ has i.i.d.~components whose common distribution satisfies the condition that 
\begin{eqnarray*}
\mathbb{E}[e^{t X_1}]  =  \frac{H(t)}{(\lambda - t)^\alpha}, \  \ \textrm{for all $t < \lambda$}
\end{eqnarray*}    
for some $\alpha > 0, \lambda > 0$ and a continuous function $H \ge 0$ such that $\lim_{t \searrow \lambda}  H(t)  \in (0, \infty)$, then the model in (\ref{ULR}) is identifiable up to a permutation of the components of $\beta_0$, i.e., in terms of \eqref{eq:plain_identifiability} the implication $\beta_0 = \beta_1$ is replaced with $\beta_{0 \, \pi(i)} = \beta_{1{i}}$ for some permutation $\pi: \{1, \ldots, d\}  \to \{1, \ldots, d\}$.


In this note, we thoroughly review the i.i.d.~case before studying the non-i.i.d.~case. We show that even when the components $X_1, \ldots, X_d$ are assumed to be independent but do not have the same distribution, a positive result asserting identifiability under broad conditions is not easy to obtain. Some simple counterexamples will be presented to highlight the challenge. On the other hand, if certain parametric assumptions are made on the distributions of the components, then one can show that the linear model in (\ref{ULR}) is identifiable. We also present a general identifiability
result for $d = 2$ based on moments, and study the apparent connection to Independent Component Analysis (ICA,  \cite{Hyvarinen2009}), a widely used tool for blind signal separation (BSS, \cite{Comon2010handbook}). Informally, ULR can be seen as a specific instance of (non-blind) source separation in which the signal to be decomposed is univariate. 
\vskip1ex
In brief, the purpose of this note is to stimulate interest in the problem at hand and equip interested reader with an overview on the current state 
of knowledge and existing obstacles towards affirmative answers, in a self-contained manner. Our hope is 
that our exposition will further propel advances in the emerging field of unlinked data analysis. 


\section{Problem Statement}\label{statement}
Consider the ULR model in \eqref{ULR}. When the goal is inference for the regression vector, a basic prerequisite is identifiability. The latter property can be formulated based on the set 
\begin{eqnarray*}
\mathcal B_0 =  \left \{ \beta \in \mathbb R^d:  Y \overset{d}{=} \beta^{\T} X + \eps  \right \}.   
\end{eqnarray*}
The following result allows us to eliminate the dependence on $\epsilon$. 

\medskip
\begin{lemma}\label{meister}
If the set of zeros of the characteristic function of $\epsilon$ does not contain any open, nonempty interval as a subset, then we have the equivalence 
$\beta_0, \beta_1  \in \mathcal{B}_0  \Longleftrightarrow  \beta^{\T}_0  X  \overset{d}{=}  \beta^{\T}_1  X$.
\end{lemma}
\noindent The condition regarding the distribution of $\eps$ in Lemma \ref{meister} is mild, and is satisfied by all discrete distributions and those commonly used in regression settings.  Examples include the Gaussian, Laplace, symmetric Gamma, t-distribution, uniform distributions, and any finite convolution thereof.


\noindent Under the sufficient condition of Lemma \ref{meister}, we have that
\begin{eqnarray}\label{B0}
\mathcal B_0 =  \left \{ \beta \in \mathbb R^d:  \beta^{\T} X \overset{d}{=} \beta^{\T}_0 X   \right \}.   
\end{eqnarray}
for any fixed element $\beta_0 \in \mathcal{B}_0$.  Depending on the size of $\mathcal B_0$ we can distinguish between the following cases:

\begin{itemize}
\item  $| \mathcal{B}_0 | = \infty$ -- non-identifiability;

\item  $| \mathcal{B}_0 | < \infty$  -- (weak) identifiability;

\item $| \mathcal{B}_0 | = 1 $  -- (strong) identifiability.
\end{itemize}

\section{The i.i.d.~case: a review}\label{IID}
Considering the case in which $X_1,\ldots,X_d$ are i.i.d., the question studied in this paper bears a close connection to the literature on so-called characterization problems \cite{Kagan1973}. In a nutshell, such problems are concerned with properties of classes of 
(typically linear) operations on a set of independent random variables that uniquely characterize the underlying distribution. For example, if $X_1, X_2 \sim \mc{N}(0,1)$ then $X_1 + X_2$ and $X_1 - X_2$ are independent; conversely, if $X_1 + X_2$ and $X_1 - X_2$ are independent, then $X_1$ and $X_2$ must be Gaussian, i.e., independence of the sum and the difference is a {\em characterizing property}. G.~P\'{o}lya \citep{polya1923herleitung} showed that if $X$ and $Y$ are i.i.d.~random variables such that their common distribution has a bounded and Riemann-integrable density on any finite interval and admits a finite second moment, then $aX + b Y$ and $c X$ (with $a, b, c >0$) have the same distribution if and only if the common distribution of $X$ and $Y$ is Gaussian. Extensions of this result due to Marcinkiewicz \cite{marcinkiewicz1939propriete} and Linnik        \cite{linnik1953linear} stated below answer the identifiability question of the present paper in the i.i.d.~case to a significant degree.


\begin{theorem}\label{M39} (Marcinkiewicz, 1939)
\par \noindent Let $\{X_n\}$ be a sequence of i.i.d.~random variables having moments of all positive orders, and let $\{a_n\}, \{b_n\}$ be (possibly infinite) sequences of constants such that 
\begin{eqnarray*}
\sum_n a_n X_n \overset{d}{=}  \sum_n b_n X_n.   
\end{eqnarray*}
Then, either  $\{\vert a_n \vert \}_n$ is a permutation of $\{\vert b_n \vert \}_n$, or the common distribution of $\{ X_n \}$ is Gaussian.
\end{theorem}
This theorem is strong in that it asserts weak identifiability as long as the 
common distribution underlying each term in the linear combination has moments of all order -- the Gaussian distribution being the notable exception.

\bigskip

To present the result shown by Linnik in 1953 (in a slightly simplified form), we largely follow the exposition in the excellent survey on {\em Identically Distributed Linear Forms and the Normal Distribution} \cite{ghurye1973identically}. Consider real numbers $a_1,\ldots, a_d, b_1, \ldots, b_d$, and define the following function
\begin{eqnarray}\label{FunctionTau}
\tau(x)=  \sum_{j=1}^d a^{2x}_j  -  \sum_{j=1}^d b^{2x}_j, \; x \ge 0.  
\end{eqnarray}

\begin{theorem}\label{Lin53}
(Linnik, 1953)

\par \noindent   Suppose $X_1,\ldots,X_d$, are i.i.d.~random variables with common distribution. Let $a_1,\ldots, a_d$, $b_1,\ldots, b_d$ be constants, and let\\[-4.5ex]
\begin{alignat*}{2}
&a = \max\{a_1,\ldots,a_d,b_1, \ldots,b_d\}, \qquad &&p = \sum_{i=1}^d  \mathds{1}_{\{ a_i =a \}}, \qquad q = \sum_{i=1}^d  \mathds{1}_{\{b_i =a \} },\\ 
 &\xi_0 = \max \{ \xi \in (0, \infty): \tau(\xi) =0 \}.
\end{alignat*}
If $p \ne q$, then in order for the following implication 
\begin{eqnarray*}
\sum_{i=1}^d a_i X_i  \stackrel{d}{=} \sum_{i=1}^d b_i X_i  \Longrightarrow  X_i \sim \mathcal{N}(\mu, \sigma^2)  \;\, \forall i \in \{1, \ldots, d \}
\end{eqnarray*}
to hold true, it is necessary and sufficient that
\begin{itemize}
\item[{\em(A)}] every positive zero of $\tau$ defined in (\ref{FunctionTau}) which is less than $\xi_0$ is an integer and a simple zero\footnote{A zero (i.e., root) of $\tau$ is called simple if the derivative of $\tau$ evaluated at this root is non-zero.};
\item[{\em(B)}] either $\xi_0$ is a simple zero and $2m -1 \le \xi_0 \le 2m $ for some integer $m$, or $\xi_0$ is a double zero\footnote{$\xi_0$ is a double zero if $\tau(\xi_0) = \tau'(\xi_0) = 0$.} and is an odd integer.    
\end{itemize}
Furthermore, if {\em (A)} and {\em (B)} are satisfied, then $\mu=0$ unless $\sum_{i=1}^d a_i =  \sum_{i=1}^d b_i$ and $\sigma^2 =0$ unless $\sum_{i=1}^d a_i^2 = \sum_{i=1}^d b^2_i$.
\end{theorem}
 Note that in the above theorem, no moment assumption is made about the common distribution of the $X_i$'s. In this regard, Linnik's result is stronger than the one by Marcikiewicz. The two results are antithetical in the sense that the former result holds universally for all
distributions but places intricate restrictions on the coefficients $\{ a_i \}$ and $\{ b_i \}$, whereas the latter has (essentially) universal validity with regard to the coefficients, but
considerably limits the class of distributions.



\begin{figure}
\begin{center}
\includegraphics[width = .45\textwidth]{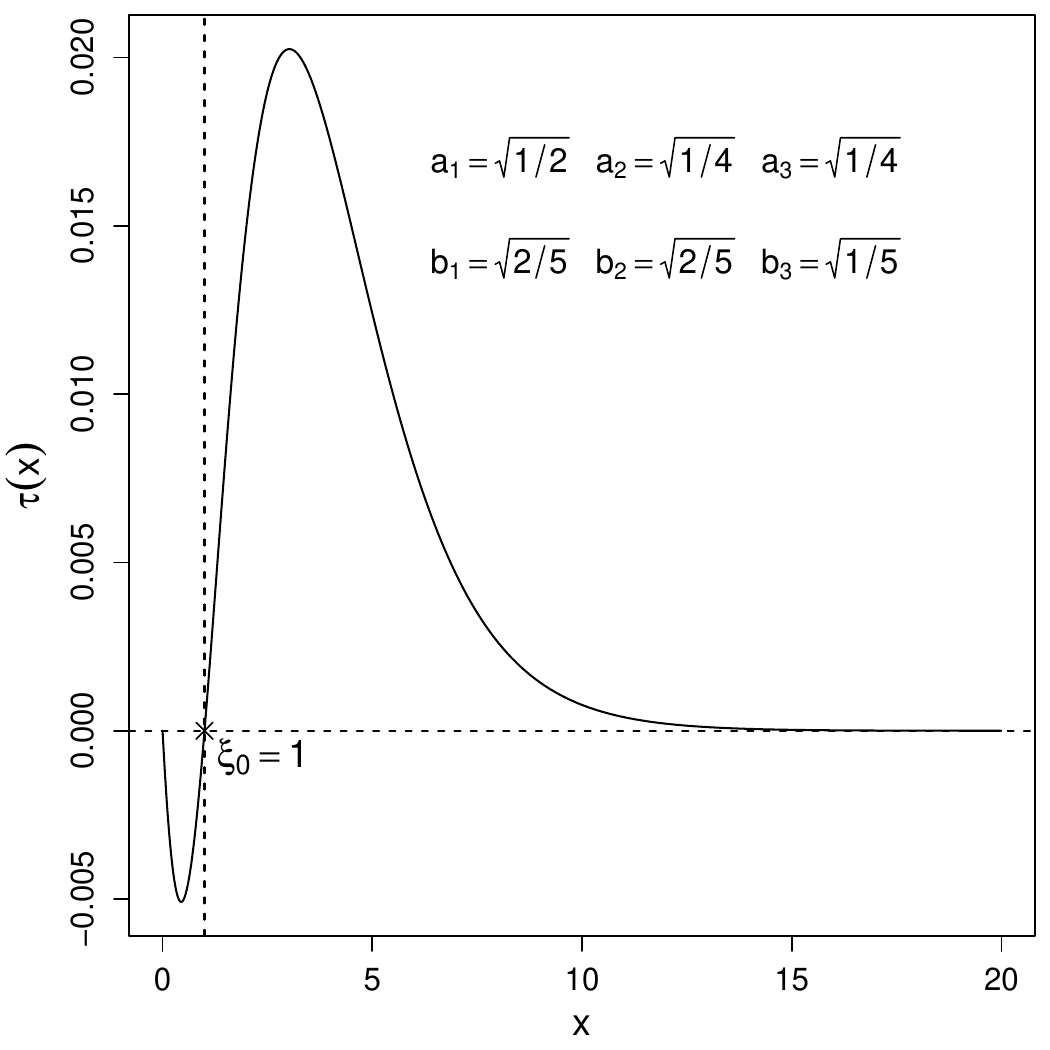} 
\hspace*{6ex} \includegraphics[width = .45\textwidth]{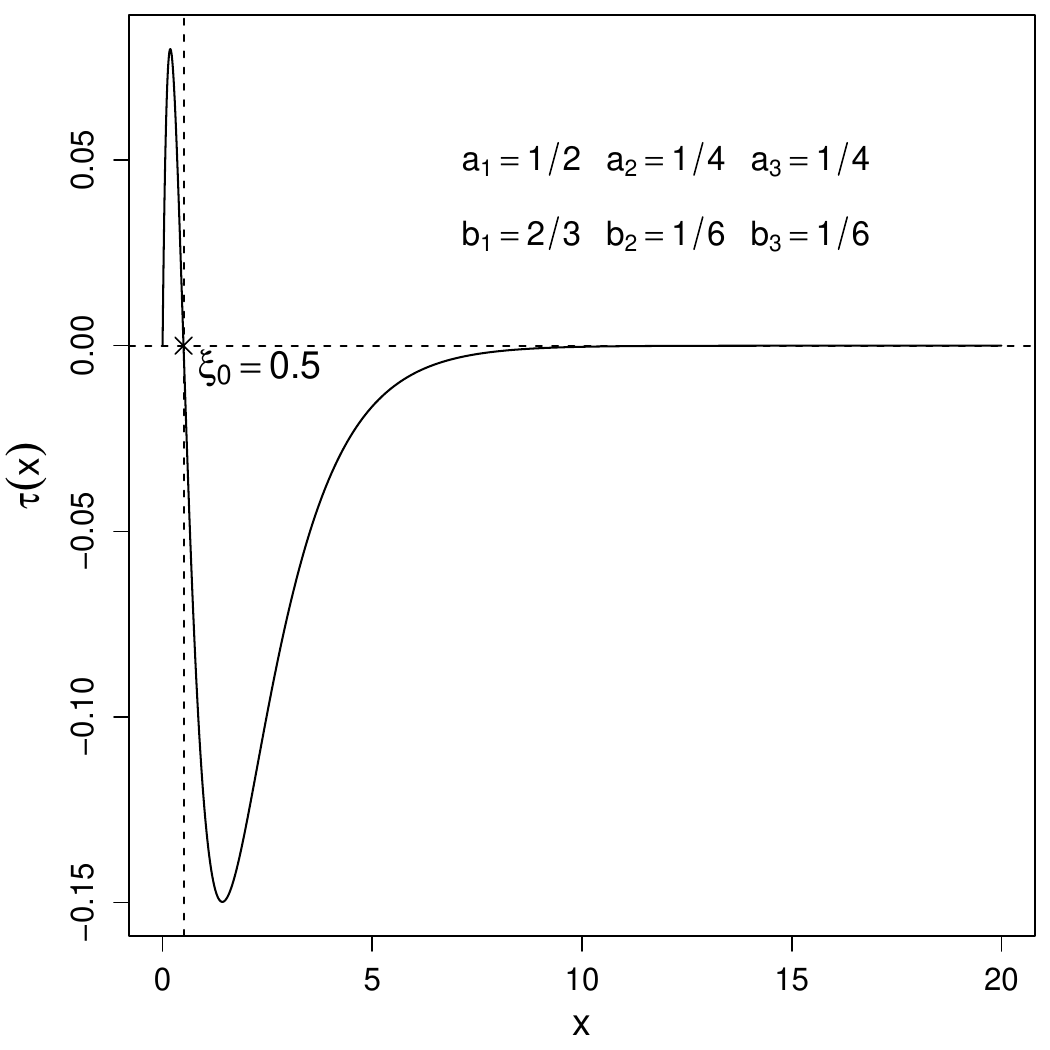} 
\end{center}
\vspace*{-4ex}
\caption{Illustration of Theorem \ref{Lin53} for $d = 3$ and two pairs of coefficients. Left: $\xi_0 = 1$ is the only root of $\tau$; since $\xi_0$ is simple and an odd integer, the conditions of the theorem are satisfied. Right: $\xi_0 = .5$ is the only root of $\tau$, but condition (B) in the theorem is {\em not met} since the integral part of $\xi_0$ is even.}\label{fig:linnik}
\end{figure}

Since investigating the zeros of the function $\tau$ is generally not straightforward  (this is certainly the case for large $d$), \cite{ghurye1973identically} provide some sufficient conditions under which  (A) and (B) hold true.  Their Lemma 1.5 gives such conditions in the following form: 
\begin{eqnarray*}
 && a^2_1  \ge \ldots  \ge a_d^2, \  b^2_1\ge \ldots \ge b^2_d, \  \textrm{and} \\
&& \sum_{i=1}^k a_i^2 \ge \sum_{i=1}^k b_i^2, \  \textrm{for $k=1, \ldots, d-1$},  \ \textrm{and} \  \ \sum_{i=1}^d a_i^2= \sum_{i=1}^d b_i^2
\end{eqnarray*}
with strict inequality for at least one $k \in \{1, \ldots, d-1\}$. Figure \ref{fig:linnik} shows an example in which the above conditions are not satisfied, but the conditions of Theorem \ref{Lin53} are seen to hold: $a_1 = \sqrt{1/2}$, $a_2 = a_3 = \sqrt{1/4}$, and 
$b_1= b_2 = \sqrt{2/5}$, $b_3 = \sqrt{1/5}$. In the second example of the figure, the sums of the coefficients are the same (and their sums of squares differ), hence the implication in Theorem \ref{fig:linnik} cannot be true -- it is well known that if the $\{ X_i \}$ are i.i.d.~standard Cauchy RV's then any linear combination of these is again a centered Cauchy random variable with scale parameter equal to the $\ell_1$-norm of the coefficients.       

\section{The non-i.i.d.~case: results and open problems }\label{NonIID}
For an integer $d \ge 2$, let us consider $d$ random variables $X_1, \ldots, X_d$ defined on the same probability space and let $X = (X_1, \ldots, X_d)^{\T}$. Suppose that there exist vectors $\beta_0$ and $\beta_1$ in $\mathbb R^d$
such that $\beta^{\T}_0  X \stackrel{d}{=}  \beta^{\T}_1  X$.  In the i.i.d.~case reviewed in $\S$\ref{IID}, it can be shown that under certain conditions, $\beta_0$ and $\beta_1$ are either equal up to some permutation and changes in signs, or the common distribution of the $X_i$'s is Gaussian. When the components $X_i, i=1, \ldots, d$, do not have the same distribution, the problem becomes much more challenging as we showcase through the results proved in the following sections.

\subsection{Scale family distributions}

The following result is a straightforward consequence of Theorem \ref{M39}. 

\begin{theorem}\label{ScaleFam}
Let  $f$ be a density function such that $f$ is not a Gaussian density, and $\int |x|^m f(x) \, dx < \infty$ for all $m > 0$.  Let $\{ f_\lambda :=\lambda f(\lambda \, \cdot), \ \lambda \in (0, \infty)\}$ be the generated scale family.  Consider $X= (X_1, \ldots, X_d)^{\T}$ to be a random vector with independent components $X_j  \sim f_{\lambda_j}, \ 1 \leq j \leq d$. Also, suppose that the noise $\epsilon$ satisfies the condition of Lemma \ref{meister}. Then, if $\beta_0 = (\beta_{01}, \ldots, \beta_{0d})^{\T}$ is a regression vector in the ULR model \eqref{ULR} it holds that
\begin{eqnarray*}
&& \Big \{ \beta \in \mathbb R^d:   \beta_j  =  \frac{\lambda_j}{\lambda_{\sigma(j)}} \beta_{0\sigma(j)} \; \textrm{for some permutation} \; \sigma:  \{1, \ldots, d\}  \mapsto \{1, \ldots, d \}  \Big\}  \\
&& \subseteq  \mathcal{B}_0 \\
&& \subseteq  \Big \{ \beta \in \mathbb R^d:  \vert \beta_j \vert =  \frac{\lambda_j}{\lambda_{\sigma(j)}} \vert \beta_{0\sigma(j)}\vert \; \textrm{for some permutation} \; \sigma:  \{1, \ldots, d\}  \mapsto \{1, \ldots, d \}  \Big\}.
\end{eqnarray*}
Moreover, in case $f$ is symmetric around $0$, then the second inclusion becomes an equality.  
\end{theorem}

\noindent We note that the assertion of the above statement continues to hold if $X$ is replaced by 
$M X$ for some invertible matrix $M \in \mathbb R^{d \times d}$.

\subsection{Elliptical symmetry}

The following result establishes that the non-identifiability of ULR for Gaussian $X$ rooting in rotational invariance extends to the entire class of spherically symmetric distributions. We note, however, that the Gaussian distribution is special within that class since 
its components are independent \cite{Fang2018}. 
\begin{theorem}\label{Spher}
Suppose that $X = (X_1, \ldots, X_d)^{\T}$ has a spherically symmetric distribution such that $X_i$ is not degenerate at $0$ for all $i = 1, \ldots, d$. Then there exists $\rho > 0$  such that
\begin{eqnarray*}
\mathcal{B}_0 =  \{ \beta \in \mathbb R^d:  \Vert \beta \Vert =  \rho \}  \invcoloneq \mathcal{S}_{\rho},
\end{eqnarray*}
the sphere of radius $\rho$.
\end{theorem}

\noindent The following corollary to Theorem \ref{Spher} yields a further extension to elliptically symmetric distributions.  
\begin{corro}\label{Ellip}
Suppose that $X = (X_1, \ldots, X_d)^{\T}$ has a elliptically symmetric distribution with parameters $\mu  =(\mu_1, \ldots, \mu_d)^{\T} \in  \R^d$ and $\Sigma \in \R^{d \times d }$ such that $\Sigma$ is positive definite. If $\mathbb P(X_i = \mu_i) <  1$ for all $i \in \{1, \ldots, d \}$, then there exists $c \in \mathbb R$ and $\rho > 0$  such that
\begin{eqnarray*}
\mathcal{B}_0 =  \{ \beta \in \mathbb R^d: \beta^{\T} \mu = c \}  \cap \{  \beta \in \mathbb R^d: \Vert \Sigma^{1/2}\beta \Vert =  \rho \}. 
\end{eqnarray*}
\end{corro}

\subsection{Convolutions}\label{subsec:convolutions}
Another source of invariance arises from convolutions. Suppose that there exists $j \in \{1,\ldots,d\}$ such that $X_j \overset{d}{=} \sum_{i \in I} \alpha_i X_i$ for $I \subseteq \{1,\ldots,d\} \setminus \{j \}$ and real numbers $\{ \alpha_i \} $, i.e., one of the components of $X$ is distributed as
a linear combination of the other components. As a result, $\mc{B}_0$ also contains 
$\beta_0 + \delta$, where $\delta_j = -\beta_{0j}$, $\delta_i = \alpha_i \beta_{0j}$, $i \in I$, and $\delta_i = 0$ for $i \notin \{j, I\}$. As a specific example, consider $X_1,X_2 \sim \textrm{Exp}(1)$ and $X_3 \sim \textrm{Gamma}(2, 1)$ with $\{ X_i \}_{i = 1}^3$ independent, so that $X_3 \overset{d}{=} X_1 + X_2$. 
Clearly, there are several other parametric families of distributions that are closed under convolution such as the Poisson, Cauchy, Binomial and Negative Binomial (when restricted to a common probability of success).



\subsection{Example of a class of distributions for which strong identifiability holds}
In the preceding subsections, we have pointed out specific obstacles in the way of identifiability. In the present subsection, we present a positive example for which strong identifiability  can be established. This example arises as a
generalization of Example 1 in \cite{monafadoua} in which $X=(X_1, X_2)^{\T}$ 
with $X_1 \sim \mathcal{N}(1,1)$ and $X_2 \sim \text{Exp}(1)$ such that 
$X_1$ and $X_2$ are independent. Using arguments based on moment-generating functions, 
strong identifiability is shown to hold in \cite{monafadoua}. In the statement below, this
example is generalized by (i) allowing $\mathcal{N}(\mu, \sigma^2)$ for any $\mu \ne 0$ and $\sigma > 0$, and (ii) allowing any finite collection of Gamma distributions satisfying 
certain constraints.  


\medskip

\begin{theorem}\label{Firstclass}
 Let $d \ge 3$. Consider the ULR model in (\ref{ULR}) with $\beta_0 = (\beta_{01},\ldots, \beta_{0d})^{\T}$ and $X = (X_1,\ldots,X_d)^{\T}$  such that 
    \begin{itemize}
        \item $X_1,\ldots,X_d $ are independent,
        \item $X_{d} \sim\mathcal{N}(\mu,\sigma^2), \mu \neq 0, \sigma > 0$,
        \item For $i = 1,\ldots,d-1$, $X_i \sim \mathrm{Gamma}(\alpha_i,\lambda_i)$ with $\alpha_i > 0$ and $\lambda_i > 0$,
        \item Whenever $\sum\limits_{i\in I}\alpha_i = \sum\limits_{j\in J}\alpha_j$ for subsets $I$ and $J$ of $\{1,\ldots,d-1\}$, then $I = J$.
    \end{itemize}
    Under these conditions, it holds that $\mathcal{B}_0 = \{\beta_0\}$.
\end{theorem}
\noindent While from a perspective of pursuing general results the setting in the above theorem may appear rather specific, it is worth pointing out that its proof requires significant effort. This may be an indication that proving general (positive) identifiability results is either intrinsically hard, or different tools are needed to establish such results.

\subsection{Identifiability for $d = 2$ using fourth moments}
Clearly, a necessary condition for $\beta$ to be contained in $\mc{B}_0$ is that 
$\E[(\beta^{\T} X)^m] = \E[(\beta_0^{\T} X)^m]$ for all positive integers $m$ such 
that $\E[|X_j|^m]$ is finite, $j=1,\ldots,d$. Introducing the moment tensor \cite{Mccullagh2018} $T^{(m)} = \E[X^{\otimes^m}] = \big(\E[X_{j_1} \ldots X_{j_m}] \big)$ with $1 \leq j_{\ell} \leq d$, 
$1 \leq \ell \leq m$, this can be expressed equivalently as the condition 
\begin{equation}\label{eq:momentensor}
\nscp{T^{(m)}}{\beta^{\otimes^{m}}} = \nscp{T^{(m)}}{\beta_0^{\otimes^{m}}}, 
\end{equation}
where here and above $^{\otimes^{m}}$ denotes the $m$-th order tensor product of a vector with itself, e.g., 
$\beta^{^{\otimes^m}} = \beta \overset{\text{1st}}{\otimes} \ldots \overset{m\text{th}}{\otimes} \beta = \big( \beta_{j_1} \cdot \ldots \cdot \beta_{j_m} \big)$, $1 \leq j_{\ell} \leq d$, $1 \leq \ell \leq m$, and 
$\nscp{\cdot}{\cdot}$ denotes the tensor inner product, i.e., for tensors $A$ and $B$ of the same dimensions $\nscp{A}{B} = \sum_{j_1 = 1}^d \ldots \sum_{j_m = 1}^d A_{j_1 \ldots j_m} B_{j_1 \ldots j_m}$. Eq.~\eqref{eq:momentensor} defines an $m$-th order polynomial equation in $\beta$. Solving such an equation is not straightforward in general. In this subsection, we study this problem for $d = 2$. We show that in this case, identifiability can be assessed via a simple, explicit condition involving $\beta_0$ and the second and fourth moments of $X = (X_1, X_2)^{\T}$. A reduction to moments necessarily leads to a loss of information and thus to a loss of sharpness in the sense that the condition in Theorem \ref{prop: nonIID4thmom} below may not be satisfied even though strong identifiability holds. At the same time, the following theorem can still be used to recover 
earlier results such as Example 1 in \cite{monafadoua}.

To simplify the statement of our result, we shall assume that 
$\E[X_1] = \E[X_2] = 0$ and that $\E[X_1^2] = \E[X_2^2] = 1$. In view of Lemma \ref{lem:centering_scaling} below, these assumptions can be made since centering and re-scaling yields a superset of $\mc{B}_0$, i.e., none of the solutions are lost in the process. Moreover, to avoid notational clutter, we write $(\alpha, \beta)$ instead of $\beta_0^{\T} = (\beta_{01}, \beta_{02})$.

\begin{lemma}\label{lem:centering_scaling}
Let $X_1$ and $X_2$ be independent random variables with means $\mu_1$ and $\mu_2$ and variances
$\sigma_1^2$ and $\sigma_2^2$. Let $\wt{X}_1 = X_1 - \mu_1$ and 
$\wt{X}_2 = X_2 - \mu_2$, $Z_1 = \wt{X}_1 / \sigma_1$ and $Z_2 = \wt{X}_2 / \sigma_2$. Then:
\begin{align*}
\mc{B}_0 &= \{(a, b) \in \R^2: X_1 a + X_2 b  \overset{d}{=}  X_1 \alpha + X_2 \beta \} \\
         &\subseteq \{(a, b) \in \R^2: \wt{X}_1 a + \wt{X}_2 b  \overset{d}{=}  \wt{X}_1 \alpha + \wt{X}_2 \beta \} \\
         &= \{(a, b) \in \R^2: Z_1 (a \sigma_1) + Z_2 (b \sigma_2)  \overset{d}{=}  Z_1 (\alpha \sigma_1) + Z_2 (\beta \sigma_2) \}.  
\end{align*}
\end{lemma}
\noindent The reverse containment $\supseteq$ is generally not true: take, e.g., 
$X_1 \sim \mc{N}(1, 1)$ and $X_2 \sim \mc{N}(2,1)$, then $(-\alpha, -\beta)$ is 
contained in the middle set in the above display, but is not contained in 
$\mc{B}_0$ unless $\alpha + 2 \beta = 0$. 


\begin{theo}\label{prop: nonIID4thmom}
 \ Let $X_1$ and $X_2$ be two independent, centered and unit variance random variables with finite fourth moments  $m_i  = \E[X_i^4], \ i =1,2$. Assume further that either $m_1 \neq 3$ or $m_2 \neq 3$.  Consider $Y_0 \overset{\text{d}}{=} \alpha X_1  +  \beta X_2$ such that $(\alpha, \beta) \ne (0,0)$, and its rescaled version    $Y_u = \alpha_u X_1 + \beta_u X_2$, where
$(\alpha_u, \beta_u)^{\T}=  (\alpha, \beta)^{\T} / \nnorm{(\alpha, \beta)^{\T}}_2$.\\ Let 
$\mc{B}_0 =   \{(a, b): a X_1 + b X_1 \overset{\text{d}}{=} Y_0 \}$, and define 
$$
c :=  \mathbb E[Y^4_u] = \alpha_u^4 m_1 + \beta_u^4 m_2 + 6 \alpha_u^2 \beta_u^2. 
$$
\begin{itemize}
    \item  If $c = 3$, then $\mc{B}_0 \subseteq \{(\sigma \alpha, \tau \beta), \; \sigma \in \{-1,1\}, \, \tau \in \{-1,1\}  \}$. 
    \item If $c \ne 3$, then $|\mc{B}_0| \leq 8$. Moreover, consider 
    \begin{equation*}
      w_1 = \frac{m_1 - 3}{c - 3},  \qquad w_2 = \frac{m_2 - 3}{c-3}.
\end{equation*}
If $\min\{w_1, w_2\} < 1$, then $\mc{B}_0 \subseteq \{(\sigma \alpha, \tau \beta), \; \sigma \in \{-1,1\}, \, \tau \in \{-1,1\}  \}$. 
\end{itemize}
\end{theo}
\noindent We note that the theorem states that the set $\mc{B}_0$ has at most eight elements unless the fourth moments of both $X_1$ and $X_2$ are equal to three, which in particular is the case if $X_1$ and $X_2$ are independent $\mc{N}(0,1)$-random variables. Under additional conditions, identifiability holds up to a change of signs of the coefficients. As demonstrated in the sequel, these conditions are easy to check given the fourth moments of $X_1$ and $X_2$ and the coefficients $(\alpha, \beta)$. If $X_1$ and $X_2$ are i.i.d., these conditions are
not met since in this case $w_1 = w_2 = w \ge 1$ (where the inequality is established in Lemma \ref{lem:prop_weights} in the appendix). This is expected since in the i.i.d.~case sign flips {\em and} permutations are sources of invariance, yielding up to $2! \cdot 2^2 = 8$ solutions, in which case the bound $|\mc{B}_0| \leq 8$ is attained.  
\vskip1ex
\noindent {\em Examples}. 
\begin{itemize}
\item[1.] (cf.~Example 1 in \cite{monafadoua}). Let $X_1 \sim \mc{N}(1,1)$ and $X_2 \sim \mathrm{Exp}(1)$. The fourth moment of the zero
          mean, unit variance random variable $X_2 - 1$ is equal to $9$, in particular, it is different from $3$. Since after centering $X_1$, we have $m_1 = 3$ and thus $w_1 = 0$, Theorem \ref{prop: nonIID4thmom}
          yields identifiability up to sign flips. Since $\E[\alpha X_1 + \beta X_2] = \alpha + \beta$, changing signs of $\alpha$, $\beta$, or both $\alpha$ and $\beta$, yields the conditions $\alpha = 0$, $\beta = 0$, or 
          $\alpha + \beta = 0$, respectively. In the first two cases, strong identifiability can be easily verified. 
          In the third case, we also must have $\E[(\alpha X_1 + \beta X_2)^3] = 0$, which, in view of $\alpha = -\beta$ yields 
          \begin{equation*}
          \E[X_2^3] - \E[X_1^3] + 3\E[X_1^2] \E[X_2] - 3\E[X_2^2] \E[X_1] = 0
          \end{equation*}
      The left hand side, however, evaluates as $2$, yielding a contradiction. This shows that either $\alpha \beta \ne 0$ and no sign flip is possible, or that $(0, \beta)$ (resp.~$(\alpha, 0)$) is the only solution. Hence, strong identifiability holds in this case. 
\item[2.] Let $X_1 \sim \mathrm{Laplace}(0, 1/\sqrt{2})$ and $X_2 \sim \mathrm{Uniform}(-\sqrt{3}, \sqrt{3})$, so that $\E[X_i] = 0$ and $\E[X_i^2] = 1$ for $i = 1, 2$. Also, we have
  $\E[X_1^4] = 6$, $\E[X_2^4] = 9/5$. It can be checked numerically that for any $(\alpha, \beta)$ such that $\alpha^2 + \beta^2 = 1$, at least one of the associated weights $w_1$, $w_2$ is negative, yielding uniqueness up to changes
  in signs of $\alpha$ and $\beta$. 
\item[3.] Let $X_1 \sim \mathrm{Laplace}(0, 1/\sqrt{2})$ and $X_2 \sim \mathrm{t}_5(0, \sqrt{3/5})$, where $\mathrm{t}_{\nu}(\mu, \sigma)$ denotes the t-distribution with $\nu$ degrees of freedom, location $\mu$ and scale $\sigma$. Numerically, it can be shown that $\min\{w_1, w_2 \} \geq 1$ 
for any $(\alpha, \beta)$, $\alpha^2 + \beta^2 = 1$ such that $|\alpha| \geq .1835$, i.e., 
Theorem \ref{prop: nonIID4thmom} generally does {\em not} guarantee identifiability up to signs. 
\end{itemize}
While in the first and second example, $\mc{B}_0$ is correctly characterized, we conjecture that 
identifiability up to signs holds true in the third example as well, but the utilization of moments up to order four alone may be insufficient to confirm this. 
\vskip1ex
\noindent {\em Extension to more than two random variables}. Theorem \ref{prop: nonIID4thmom} can be used for cases in which $d > 2$, but its application becomes more tedious as $d$ increases. For 
$d = 3$, Theorem \ref{prop: nonIID4thmom} can be applied for all partitions of the index
set $\{1,2,3\}$ into non-empty subsets, i.e., $\{\{1,2\}, \{3\}\}$, $\{\{1,3\}, \{2\}\}$,
and $\{\{2,3\}, \{1\}\}$. For each of these three partitions, we investigate identifiability for i) the two-element subset,
ii) the linear combination (associated with its respective regression coefficients) of the $X$'s in the two-element subset and
the remaining $X$ in the one-element subset. As a result, for each partition Theorem \ref{prop: nonIID4thmom} is invoked
twice, yielding six instances in total. If the fourth moments are not both equal to three in each 
instance, Theorem \ref{prop: nonIID4thmom} implies that the number of solutions is upper bounded 
by $3! \cdot 2^3  = 48$. This number coincides with the maximum possible number of solutions in the i.i.d.~case, apart from the exceptions according to Theorems \ref{M39} and \ref{Lin53}. 

Albeit tedious, the approach can be extended to arbitrary $d$: the set $\{1,\ldots,d\}$ can be partitioned into
subsets of size 2 and $d-2$, respectively. We then recursively divide the subsets of size $d-2$ until the size of the subset reaches 2 (if $d$ is even) or one (if $d$ is odd).


\subsection{Connection to Independent Component Analysis (ICA)}\label{subsec:ICA}
Consider the following modification of the ULR problem \eqref{ULR} in which 
\begin{equation}\label{eq:ICA_intro}
Y^{(\ell)} \overset{d}{=} X^{\T} \beta_0^{(\ell)} + \eps^{(\ell)}, \quad \ell=1,\ldots,m, \;\; \Longleftrightarrow \, \M{Y} \overset{d}{=} \M{B}_0 X + \bm{\eps} 
\end{equation}
with $\{ \beta_0^{(\ell)} \} \subset \R^{d}$ forming the $m$ rows of the matrix $\M{B}_0$, and 
$\bm{\eps} = (\eps^{(\ell)})_{\ell = 1}^m$ is a zero-mean random vector independent of $X$. 
The setup \eqref{eq:ICA_intro} is simply that of unlinked linear regression with multiple 
response variables $\M{Y} = (Y^{(\ell)})_{\ell = 1}^m$. If $d = m$ and the components of $X$ are assumed to be independent while $\M{B}_0$ is assumed to be invertible, the setup becomes that of classical {\em independent component analysis} (ICA), a widely used tool in signal processing and machine learning with a voluminous body of literature \cite[e.g.,][]{Comon1994, Hyvarinen2000, Hyvarinen2009}. Different from ULR, ICA aims to recover $\M{B}_0$ {\em and} the ``sources" $(X_i)_{i = 1}^d$ given $\M{Y}$. Therefore, in the ICA literature the question of identifiability is studied with respect to $\M{B}_0$ and $X$ jointly. Independence and non-Gaussianity of the components arise as the key requirements in that literature even though it has been shown recently that the former requirement can be relaxed \cite{mesters2024}. The ULR setup can be seen as a specific case of that under consideration in {\em overcomplete ICA} in which $m < d$. The following result concerning identifiability of $\M{B}_0$ in overcomplete ICA can be found in Theorem 1 in \cite{Eriksson2004}, which in turn references a result in \cite{Kagan1973}. 
\begin{theorem}\label{ICA_overcomplete-thm} {\em \cite[][Theorem 10.3.5]{Kagan1973}} 
Let $\M{B}_0$ and $\M{B}$ be $m$-by-$d$ matrices, $m \geq 2$, so that no pair of columns are
linearly dependent, i.e., for all $j \neq k$, there does not exist a $\lambda \in \mathbb{R}$ such that 
$\M{B}_{0j} = \lambda \M{B}_{0k}$ or $\M{B}_{j} = \lambda \M{B}_{k}$. Let further $X$ and $Z$ be $d$-dimensional random vectors with non-Gaussian, independent, and unit variance components. Then $\M{B}_0  X \stackrel{d}{=} \M{B} Z$ implies that 
$\M{B} = \M{B}_0 \M{S} \bm{\Pi}  $ for a permutation matrix $\bm{\Pi}$ and a diagonal matrix $\M{S}$ whose diagonal entries take values in $\{-1,1\}$. 
\end{theorem}
\noindent Note that in the above theorem, no statement is made regarding the identifiability 
of $X$. Letting $X = Z$ the theorem implies weak identifiability of the multi-response variant of ULR given in \eqref{eq:ICA_intro}, i.e., 
the set $\mathscr{B}_0 = \{\M{B} \in \R^{m \times d}: \; \M{B} X \overset{d} =  \M{B}_0 \M{X} \}$
contains at most $\M{B}_0$ and matrices obtained through permuting or flipping signs of its columns. 

In Theorem \ref{ICA_overcomplete-thm} {\em all} source variables are required to be non-Gaussian, which may appear overly strict since results on the identifiability in classical ICA ($d = m$ with $\M{B}_0$ invertible) allow (at most) one Gaussian component in $X$. Very recently, efforts have been made to extend Theorem \ref{ICA_overcomplete-thm} to accommodate a single Gaussian source \cite{Wang2024identifiability}. Moreover, it is noteworthy that identifiability of $\M{B}_0$ requires that none of its columns is a scalar multiple of another column, which is significantly weaker than the linear independence of the entire set of columns of $\M{B}_0$. The former assumption is minimal since without it, one may drop an entire column and absorb it into another column: if the $(X_i)_{i = 1}^d$ are i.i.d.~and $\M{B}_{01} = \lambda\M{B}_{02}$, we can construct the matrix
$\M{B}$ with $\M{B}_{1} = (1 + \lambda) \M{B}_{01}$, $\M{B}_{2} = \M{0}$ and $\M{B}_{j} = \M{B}_{0j}$, $i=3,\ldots,d$, such that $\M{B}_0 X \stackrel{d}{=} \M{B} \M{X}$. 

In conclusion, the multi-response variant of ULR admits a general identifiability result that currently does not seem to be in reach for (ordinary) ULR. Intuitively, a crucial distinction between the two setups is that once $m \geq 2$ additional constraints arise for possible alternative solutions $\M{B}$ to ensure that not only the resulting marginal distributions but also all joint distributions of $\M{B} X$ agree with those of $\M{B}_0 X$. To illustrate this point, consider independent random variables $X_1, X_2, X_3 \sim \mathrm{Exp}(1)$ and $X_4 \sim \mathrm{Gamma}(2, 1)$ and 
\begin{equation*}
\underbrace{\begin{bmatrix}
    1 & 1 & 0 & 0\\
    0 & 1  & 1 & 0\\
    0 & 0 & 1 & 1 
    \end{bmatrix}}_{\M{B}_0} \begin{bmatrix}
    X_1 \\
    X_2 \\
    X_3 \\
    X_4
    \end{bmatrix}
\end{equation*}
Note that the columns of $\M{B}_0$ are pairwise linearly independent, and thus Theorem \ref{ICA_overcomplete-thm} applies. On the other hand, none of the individuals rows of $\M{B}_0$ are (weakly) identifiable  when considering the components of $\M{B}_0 X$ in isolation, in view of the invariance under convolutions as discussed in $\S$\ref{subsec:convolutions}.

\section{Discussion}\label{discussion}
In this note, we have discussed the problem of identifiability in unlinked linear regression, 
an emerging paradigm for addressing challenging data integration problems. Our study has led 
us to the re-discovery of classical results dating back more than half a century that
are important milestones in the literature on {\em Characterization Problems}. Interestingly, this literature has played a crucial role in proving identifiability of Independent Component Analysis, which has become a mainstream data analysis tool over the years. ICA and more generally, linear
structural equation models, are mature yet vivid fields of research. We hypothesize 
that more progress on the question raised herein can be made by borrowing techniques developed in those areas. We believe that approaches relying on moment and cumulant tensors as used recently in \cite{mesters2024} to push the boundary regarding identifiability of ICA-style decompositions constitute a promosing path forward. 

We conclude this discussion with an informal conjecture that would assert weak identifiability of unlinked regression for independent $X$'s in a broad sense. 
\vskip1ex
{\em If $X$ has independent, unit variance components, out of which at most one is Gaussian, and 
the representation $Y_0 \overset{d}{=} X^{\T} \beta_0$ is minimal, the regression parameter $\beta_0$ is identifiable up to sign flips and permutations of  its components.}
\vskip1ex
Here, minimality means that there does not exits a sub-vector of $X$ such that $Y_0$ can
be expressed as a linear combination of that sub-vector. This assumption is necessary to 
address invariances resulting from convolutions. The assumption of unit variance rules 
out the invariances arising in the context of scale families, and also rules out stable
distributions \cite{Zolotarev1986} other than the Gaussian. 

%

\bibliographystyle{plain}
\bibliography{identifiability/NonIID_reformatted}

\vskip6ex

\subsection*{Appendix: Proofs}

\subsubsection*{Proof of Lemma \ref{meister}}
For some random variable $Z$, let $\Phi_Z$  denote its characteristic function; i.e., $\Phi_Z(t) =\mathbb E[\exp(i t Z )], \ t \in \mathbb R$. Since $X$ and $\epsilon$ are assumed to be independent, it holds that
\begin{eqnarray*}
\beta_0, \beta_1  \in \mathcal{B}_0  \Longleftrightarrow  \Phi_{\beta^{\T}_0  X}  \Phi_\epsilon =  \Phi_{\beta^{\T}_1  X}  \Phi_\epsilon.
\end{eqnarray*}
It follows from the comments below display (2.23) on page 25 in \cite{meister2009}  that a sufficient condition for 
$$
\Phi_{\beta^{\T}_0  X}  \Phi_\epsilon =  \Phi_{\beta^{\T}_1  X}   \Longleftrightarrow \Phi_{\beta^{\T}_0  X} = \Phi_{\beta^{\T}_  X}
$$
to hold is that the set of zeros of $\Phi_\epsilon$ does not contain any open, nonempty interval as a subset. This finishes the proof. \hfill \(\square\)

\subsubsection*{Proof of Theorem \ref{ScaleFam}}
Fix $\beta_0 \in \mathcal{B}_0$ and let $\beta$ be such that there exists a permutation $\sigma: \{1, \ldots, d \} \mapsto  \{1, \ldots, d \}$ satisfying
$$
\beta_j =  \frac{\lambda_j}{\lambda_{\sigma(j)}} \beta_{0\sigma(j)}
$$
for $j= 1, \ldots, d$.  Then,
\begin{eqnarray*}
\beta^{\T} X  & =  &   \frac{\beta_1}{\lambda_1} (\lambda_1 X_1) + \ldots + \frac{\beta_d}{\lambda_d} (\lambda_d X_d)   \\
& = &  \frac{\beta_{0\sigma(1)}}{\lambda_{\sigma(1)}} (\lambda_1 X_1) + \ldots + \frac{\beta_{0\sigma(d)}}{\lambda_{\sigma(d)}} (\lambda_d X_d) \\
& \overset{d}{=} &   \frac{\beta_{0\sigma(1)}}{\lambda_{\sigma(1)}} (\lambda_{\sigma(1)} X_{\sigma(1)}) + \ldots + \frac{\beta_{0\sigma(d)}}{\lambda_{\sigma(d)}} (\lambda_{\sigma(d)} X_{\sigma(d)})
\end{eqnarray*}
using the fact that $\lambda_1 X_1, \ldots, \lambda_d X_d$ are i.i.d. random variables.  Hence, 
$$
\beta^{\T} X   \overset{d}{=} \beta_{0\sigma(1)} X_{\sigma(1)} + \ldots + \beta_{0\sigma(d)} X_{\sigma(d)}  =   \beta_0^{\T} X.
$$
By Lemma \ref{meister}, this implies that $\beta \in \mathcal B_0$, and hence the first inclusion is proved. To show the second inclusion,  consider $\beta  \in \mathcal{B}_0$. Then, $\beta^{\T}_0  X   \overset{d}{=} \beta^{\T}  X$, or equivalently
\begin{eqnarray*}
  \frac{\beta_{01}}{\lambda_1} (\lambda_1 X_1) + \ldots +  \frac{\beta_{0d}}{\lambda_d} (\lambda_d X_d)  \overset{D}{=}   \frac{\beta_{1}}{\lambda_1} (\lambda_1 X_1) + \ldots +  \frac{\beta_{d}}{\lambda_d} (\lambda_d X_d).
\end{eqnarray*}
Since $\lambda_1 X_1, \ldots, \lambda_d X_d $ are i.i.d. $\sim f$, and $f$ is a not the density of a Gaussian, Theorem \ref{M39} implies that there must exist a permutation $\sigma : \{1, \ldots, d\}  \to \{1, \ldots, d \}$ such that
\begin{eqnarray}\label{perm}
\frac{\vert \beta_j \vert}{\lambda_j}  =  \frac{\vert \beta_{0\sigma(j)}\vert}{\lambda_{\sigma(j)}}
\end{eqnarray}
for $j =1, \ldots, d$. This establishes the claimed inclusion. 

Suppose now that $f$ is even and let $\beta \in \mathbb R^d$ such that there exists a permutation $\sigma$ for which the equality in (\ref{perm}) is satisfied.  Using the symmetry of the $f$ around $0$ and independence of the $X_i'$s it follows that  
\begin{align*}
&\frac{\beta_{0\sigma(1)}}{\lambda_{\sigma(1)}} (\lambda_{\sigma(1)} X_{\sigma(1)}) + \ldots +  \frac{\beta_{0\sigma(d)}}{\lambda_{\sigma(d)}} (\lambda_{\sigma(d)} X_{\sigma(d)}) \\
&\overset{d}{=}    \frac{\vert \beta_{0\sigma(1)} \vert }{\lambda_{\sigma(1)}} (\lambda_{\sigma(1)} X_{\sigma(1)}) + \ldots +  \frac{\vert \beta_{0\sigma(d)}\vert}{\lambda_{\sigma(d)}} (\lambda_{\sigma(d)} X_{\sigma(d)})  \\
&=   \frac{\vert \beta_{1} \vert }{\lambda_{1}} (\lambda_{\sigma(1)} X_{\sigma(1)}) + \ldots +  \frac{\vert \beta_{d}\vert}{\lambda_d} (\lambda_{\sigma(d)} X_{\sigma(d)})  \\
&\overset{d}{=}  \frac{\vert \beta_{1} \vert }{\lambda_{1}} (\lambda_1 X_{1}) + \ldots +  \frac{\vert \beta_{d}\vert}{\lambda_d} (\lambda_{d} X_{d})  \\
&\overset{d}{=}   \frac{\beta_{1} }{\lambda_{1}} (\lambda_1 X_{1}) + \ldots +  \frac{\beta_{d}}{\lambda_d} (\lambda_{d} X_{d})  = \beta^{\T} X.
\end{align*}
Since
\begin{eqnarray*}
\frac{\beta_{0\sigma(1)}}{\lambda_{\sigma(1)}} (\lambda_{\sigma(1)} X_1) + \ldots +  \frac{\beta_{0\sigma(d)}}{\lambda_{\sigma(d)}} (\lambda_{\sigma(d)} X_d)  = \beta^{\T}_0  X
\end{eqnarray*}
we conclude that $\beta^{\T}_0 X \overset{d}{=}\beta^{\T} X $, and hence $\beta \in \mathcal{B}_0$. This finishes the proof. \hfill \(\square\)
$$
$$

\subsubsection*{Proof of Theorem \ref{Spher}}
$X$ is spherically symmetric if and only if for any orthogonal matrix $P  \in \mathbb R^{d \times d}$ it holds that $P X \stackrel{d}{=}  X$.  Let $\beta_0 \in \mathbb R^d$ be a regression vector in the unlinked linear regression model as in (\ref{ULR}) and let $\Vert \beta_0 \Vert := \rho$. It is well-known that $\beta \in \mathcal{S}_{\rho}$ if and only if there exists an orthogonal matrix $P  \in \mathbb R^{d \times d}$ such that $\beta =  P\beta_0$.  Then,  $\beta^{\T} X =  \beta^{\T}_0  PX \stackrel{d}{=}  \beta^{\T}_0 X$.  This means that $\mathcal{S}_{\rho}  \subset \mathcal{B}_0$. Let $\beta \in \mathbb R^d$  such that $\beta^{\T} X  \stackrel{d}{=}  \beta^{\T}_0  X$ and suppose that $\Vert \beta  \Vert  \neq \Vert \beta_0 \Vert$. Let $e_1$ be the first canonical basis vector.  Pick orthogonal matrices $ P$ and $Q$ such that
\begin{eqnarray*}
P \beta_0 =  \Vert \beta_0  \Vert e_1, \ \  \textrm{and} \  \  Q \beta =  \Vert \beta \Vert e_1.
\end{eqnarray*}
Then, $\Vert \beta_0 \Vert X_1 \stackrel{d}{=} \Vert \beta \Vert X_1$ which can only be true if $X_1$ is degenerate at $0$.  This implies that $\Vert \beta  \Vert  = \Vert \beta_0 \Vert$ and therefore $\mathcal B_0 \subset \mathcal{S}_\rho$. \hfill \(\square\)

\subsection*{Proof of Corollary \ref{Ellip}}
It follows from \cite[Theorem 1]{cambanis1981theory} that $X$ is elliptically symmetric with parameters $\mu$ and $\Sigma$ if and only if there exists $Y$ which is spherically symmetric such that
$$
X =  \mu + \Sigma^{1/2} Y.
$$
Let $\beta_0$ and $\beta_1$ in $\mathbb R^d$. Then, $\beta_0$ and $\beta_1$ belong to $\mathcal B_0$ if and only if 
\begin{eqnarray}\label{eq}
\beta_0^{\T}  \mu + \beta_0^{\T} \Sigma^{1/2} Y \stackrel{d}{=} \beta_1^{\T}  \mu + \beta_1^{\T} \Sigma^{1/2} Y.
\end{eqnarray}
Since the covariance of $Y$ is the identity, the characteristic function of $Y$ is given by
$$
\Phi_Y(u)  =  \phi(\Vert u \Vert^2), \ u \in \mathbb R^d
$$
for some real function $\phi$ which admits the representation (1) given in \cite[][p.~369]{cambanis1981theory}. Now, the equality in (\ref{eq}) is equivalent to 
\begin{eqnarray*}
 e^{it \beta_0^{\T}  \mu}  \phi(t^2 \Vert  \Sigma^{1/2} \beta_0  \Vert^2) =  e^{it \beta_1^{\T}  \mu}  \phi(t^2 \Vert  \Sigma^{1/2} \beta_1  \Vert^2)  
\end{eqnarray*}
for all $t \in \mathbb R$. Since $\phi$ is a real function, the latter means that $e^{it \beta_0^{\T}  \mu} = e^{it \beta_1^{\T}  \mu}$ for all $t  \in \mathbb R$ and hence $\beta_0^{\T}  \mu = 
  \beta_1^{\T}  \mu$. Replacing in (\ref{eq}) yields $\beta_0^{\T} \Sigma^{1/2} Y \stackrel{d}{=}  \beta_1^{\T} \Sigma^{1/2} Y$. Since
  $\mathbb P(X_j = \mu_j) < 1$ for all $j \in \{1, \ldots, d \}$ implies that none of the components $Y_j$ is degenerate at $0$, Theorem \ref{Spher} implies that $\Vert \Sigma^{1/2} \beta_0 \Vert = \Vert \Sigma^{1/2} \beta_1 \Vert$, which completes the proof.  \hfill \(\square\)

\subsubsection*{Proof of Theorem \ref{Firstclass}}
Assume $ \lvert\mathcal{B}_0\rvert > 1$. Take $\beta \neq \tilde{\beta}\in\mathcal{B}_0$ such that $\beta = (\beta_1,\ldots,\beta_d)^{\T}$ and $\tilde{\beta} = (\tilde{\beta}_1,\ldots,\tilde{\beta}_d)^{\T}$. Writing $M(t)$ for the respective moment-generating functions, we get
\[\begin{split}
    \beta^{\T}X \stackrel{d}{=}\tilde{\beta}^{\T}X & \Longleftrightarrow M_{\beta^{\T}X}(t) = M_{\tilde{\beta}^{\T}X}(t) \text{ for all }t\in\R\backslash \{0\} \\& \Longleftrightarrow \E[e^{t(\beta_1X_1 + \ldots + \beta_dX_d)}] = \E[e^{t(\tilde{\beta}_1X_1 + \ldots + \tilde{\beta}_dX_d)}] \text{ for all }t\in\R\backslash \{0\} \\& \implies M_{X_1}(\beta_1t)\ldots M_{X_d}(\beta_dt) = M_{X_1}(\tilde{\beta}_1t)\ldots M_{X_d}(\tilde{\beta}_dt) \text{ for all }t\in\R\backslash \{0\}
\end{split}\]
(using independence), whenever both sides are defined.  We arrive at the following equality which holds for all $t\in\R\backslash\{0\}$ such that both sides are defined:
\begin{eqnarray}\label{eq3.8}
    \frac{M_{X_{d}}(\beta_{d}t)}{M_{X_{d}}(\tilde{\beta}_{d}t)}\cdot M_{X_1}(\beta_1t)\ldots M_{X_{d-1}}(\beta_{d-1}t) = M_{X_1}(\tilde{\beta}_1t)\ldots M_{X_{d-2}}(\tilde{\beta}_{d-1}t)
\end{eqnarray}
where
\begin{eqnarray*}
\frac{M_{X_{d}}(\beta_{d}t)}{M_{X_{d}}(\tilde{\beta}_{d}t)} = e^{\mu(\beta_{d} - \tilde{\beta}_{d})t + \frac{\sigma^2}{2}(\beta_{d}^2 - \tilde{\beta}_{d}^2)t^2}.
\end{eqnarray*}
For $Y \sim \mathrm{Gamma}(\alpha,\lambda), \alpha > 0, \lambda > 0$ and $\gamma\in\R\backslash\{0\}$, it holds that 
\[M_Y(\gamma t) = \begin{cases}
\frac{\lambda^\alpha}{(\lambda - \lvert\gamma\rvert t)^\alpha}, t < \frac{\lambda}{\lvert\gamma\rvert} & \text{if } \gamma > 0\\
\frac{\lambda^\alpha}{(\lambda + \lvert\gamma\rvert t)^\alpha}, t > -\frac{\lambda}{\lvert\gamma\rvert} & \text{if } \gamma < 0
\end{cases}\]

Define
\[z = \lvert\{i\in\{1,\ldots,d-1\} : \beta_i = 0\}\rvert\]
and
\[\tilde{z} = \lvert\{i\in\{1,\ldots,d-1\} : \tilde{\beta}_i = 0\}\rvert,\]
the number of zero-components of $\beta$ and $\tilde{\beta}$ respectively among the first $d-1$ dimensions.  Furthermore, let 
\[I_+ = \{i\in\{1,\ldots,d-1\} : \beta_i > 0\}\]
as well as
\[I_- = \{i\in\{1,\ldots,d-1\} : \beta_i < 0\}\]
Analogously, we define
\[\tilde{I}_+ = \{i\in\{1,\ldots,d-1\} : \tilde{\beta}_i > 0\}\]
as well as
\[\tilde{I}_- = \{i\in\{1,\ldots,d-1\} : \tilde{\beta}_i < 0\}.\]
The first $d-1$ dimensions of the vectors $\beta$ and $\tilde{\beta}$ have $d-1-z$ and $d-1-\tilde{z}$ non-zero components respectively. We define the numbers of these components as
\[c_+ = \lvert I_+\rvert\]
and
\[c_- = \lvert I_-\rvert\]
for the vector $\beta$ as well as
\[\tilde{c}_+ = \lvert \tilde{I}_+\rvert\]
and
\[\tilde{c}_- = \lvert \tilde{I}_-\rvert\]
for the vector $\tilde{\beta}$. Defining  $\varphi(t) = e^{\mu(\beta_{d} - \tilde{\beta}_{d})t + \frac{\sigma^2}{2}(\beta_{d}^2 - \tilde{\beta}_{d}^2)t^2}$,  equation \eqref{eq3.8} can be written as
\begin{eqnarray}\label{eq3.9}
    \varphi(t)\cdot\prod\limits_{i\in I_+}\frac{\lambda_i^{\alpha_i}}{(\lambda_i - \lvert\beta_i\rvert t)^{\alpha_i}}\cdot\prod\limits_{j\in I_-}\frac{\lambda_j^{\alpha_j}}{(\lambda_j + \lvert\beta_j\rvert t)^{\alpha_j}} = \prod\limits_{r\in \tilde{I}_+}\frac{\lambda_r^{\alpha_r}}{(\lambda_r - \lvert\tilde{\beta}_r\rvert t)^{\alpha_r}}\cdot\prod\limits_{s\in \tilde{I}_-}\frac{\lambda_s^{\alpha_s}}{(\lambda_s + \lvert\tilde{\beta}_s\rvert t)^{\alpha_s}}
\end{eqnarray}
for $t\neq 0$, $t\in \Big(\max\limits_{j\in I_-}\big(-\frac{\lambda_j}{\lvert\beta_j\rvert}\big),\min\limits_{i\in I_+}\big(\frac{\lambda_i}{\lvert\beta_i\rvert}\big)\Big)$ as well as $t\in \Big(\max\limits_{s\in \tilde{I}_-}\big(-\frac{\lambda_s}{\lvert\tilde{\beta}_s\rvert}\big),\min\limits_{r\in \tilde{I}_+}\big(\frac{\lambda_r}{\lvert\tilde{\beta}_r\rvert}\big)\Big)$.
Define
\[R_i^+ = \frac{\lambda_i}{\lvert\beta_i\rvert},\quad i\in I_+\]
and analogously
\[R_i^- = \frac{\lambda_i}{\lvert\beta_i\rvert},\quad i\in I_-\]
as well as
\[\tilde{R}_i^+ = \frac{\lambda_i}{\lvert\tilde{\beta}_i\rvert},\quad i\in \tilde{I}_+\]
and
\[\tilde{R}_i^- = \frac{\lambda_i}{\lvert\tilde{\beta}_i\rvert},\quad i\in \tilde{I}_-.\]
Writing $R_{(i)}^+$ for the $i$-th smallest value of $R_1^+,\ldots,R_{c_+}^+$, and analogously for $R_{(i)}^-, \tilde{R}_{(i)}^+$ and $\tilde{R}_{(i)}^-$, equation \eqref{eq3.9} can be written as
\begin{equation}\label{eq3.10}
    \varphi(t)\cdot\prod\limits_{i\in I_+}\frac{1}{\big(1 - \frac{t}{R_i^+}\big)^{\alpha_i}}\cdot\prod\limits_{j\in I_-}\frac{1}{\big(1 + \frac{t}{R_j^-}\big)^{\alpha_j}} = \prod\limits_{r\in \tilde{I}_+}\frac{1}{\big(1 - \frac{t}{\tilde{R}_r^+}\big)^{\alpha_r}}\cdot\prod\limits_{s\in \tilde{I}_-}\frac{1}{\big(1 + \frac{t}{\tilde{R}_s^-}\big)^{\alpha_s}}
\end{equation}
for $t\neq 0$, $t\in \big(-R_{(1)}^-, R_{(1)}^+\big)$ as well as $t\in \big(-\tilde{R}_{(1)}^-, \tilde{R}_{(1)}^+\big)$. Note that we made use of the relation $\max\limits_{j\in I_-}\big(-\frac{\lambda_j}{\lvert\beta_j\rvert}\big) = -\min\limits_{j\in I_-}\big(\frac{\lambda_j}{\lvert\beta_j\rvert}\big) = -R_{(1)}^-$ and analogously for $\max\limits_{s\in \tilde{I}_-}\big(-\frac{\lambda_s}{\lvert\tilde{\beta}_s\rvert}\big) = -\tilde{R}_{(1)}^-$.
Assume now that $R_{(1)}^+ \neq \tilde{R}_{(1)}^+$. We assume that $R_{(1)}^+ < \tilde{R}_{(1)}^+$. Then we obtain
\[\lim\limits_{t\nearrow R_{(1)}^+}\varphi(t)\cdot\prod\limits_{i\in I_+}\frac{1}{\big(1 - \frac{t}{R_i^+}\big)^{\alpha_i}}\cdot\prod\limits_{j\in I_-}\frac{1}{\big(1 + \frac{t}{R_j^-}\big)^{\alpha_j}} = +\infty\]
and
\[\lim\limits_{t\nearrow R_{(1)}^+}\prod\limits_{r\in \tilde{I}_+}\frac{1}{\big(1 - \frac{t}{\tilde{R}_r^+}\big)^{\alpha_r}}\cdot\prod\limits_{s\in \tilde{I}_-}\frac{1}{\big(1 + \frac{t}{\tilde{R}_s^-}\big)^{\alpha_s}} < +\infty\]
which is a contradiction to \eqref{eq3.10}. The same reasoning shows that we get a contradiction if $R_{(1)}^+ > \tilde{R}_{(1)}^+$ and we therefore have $R_{(1)}^+ = \tilde{R}_{(1)}^+$. A symmetric argument shows that $R_{(1)}^- = \tilde{R}_{(1)}^-$. We can now conclude by induction that $c_+ = \tilde{c}_+$ and $R_{(i)}^+ = \tilde{R}_{(i)}^+$ for all $i = 1,\ldots, c_+$ as well as $c_- = \tilde{c}_-$ and $R_{(i)}^- = \tilde{R}_{(i)}^-$ for all $i = 1,\ldots, c_-$. Let
\[v_{(1)}^+ < \ldots < v_{(d_+)}^+\]
denote the distinct and ordered values of $\{R_{(1)}^+,\ldots,R_{(c_+)}^+\}$. Note that by our previous considerations, $v_{(1)}^+ < \ldots < v_{(d_+)}^+$ also denote the distinct and ordered values of $\{\tilde{R}_{(1)}^+,\ldots,\tilde{R}_{(c_+)}^+\}$. Define also
\[I_j^+ = \{i\in I_+ : R_i^+ = v_{(j)}^+\}\]
and
\[\tilde{I}_j^+ = \{i\in \tilde{I}_+ : \tilde{R}_i^+ = v_{(j)}^+\}\]
for all $j\in\{1,\ldots,d_+\}$ and note that $I_+ = \bigcup\limits_{j = 1}^{d_+}I_j^+ = \bigcup\limits_{j = 1}^{d_+}\tilde{I}_j^+$ form partitions of $I_+$. Analogously, let
\[v_{(1)}^- < \ldots < v_{(d_-)}^-\]
denote the distinct and ordered values of $\{R_{(1)}^-,\ldots,R_{(c_-)}^-\}$ and define
\[I_j^- = \{i\in I_- : R_i^- = v_{(j)}^-\}\]
as well as
\[\tilde{I}_j^- = \{i\in \tilde{I}_- : \tilde{R}_i^- = v_{(j)}^-\}\]
for all $j\in\{1,\ldots,d_-\}$. Equation \eqref{eq3.10} then becomes
\begin{multline*}
\varphi(t)\cdot\Bigg(\frac{1}{1 - \frac{t}{v_{(1)}^+}}\Bigg)^{\sum\limits_{i\in I_1^+}\alpha_i}\cdots\Bigg(\frac{1}{1 - \frac{t}{v_{(d_+)}^+}}\Bigg)^{\sum\limits_{i\in I_{d_+}^+}\alpha_i}\cdot\Bigg(\frac{1}{1 + \frac{t}{v_{(1)}^-}}\Bigg)^{\sum\limits_{j\in I_1^-}\alpha_j}\cdots\Bigg(\frac{1}{1 + \frac{t}{v_{(d_-)}^-}}\Bigg)^{\sum\limits_{j\in I_{d_-}^-}\alpha_j}\\ = \Bigg(\frac{1}{1 - \frac{t}{v_{(1)}^+}}\Bigg)^{\sum\limits_{r\in \tilde{I}_1^+}\alpha_r}\cdots\Bigg(\frac{1}{1 - \frac{t}{v_{(d_+)}^+}}\Bigg)^{\sum\limits_{r\in \tilde{I}_{d_+}^+}\alpha_r}\cdot\Bigg(\frac{1}{1 + \frac{t}{v_{(1)}^-}}\Bigg)^{\sum\limits_{s\in \tilde{I}_1^-}\alpha_s}\cdots\Bigg(\frac{1}{1 + \frac{t}{v_{(d_-)}^-}}\Bigg)^{\sum\limits_{s\in \tilde{I}_{d_-}^-}\alpha_s}
\end{multline*}
for $t\neq 0$ and $t\in\big(-v_{(1)}^-, v_{(1)}^+\big)$. Bringing all the factors to the left-hand side gives us
\begin{multline}\label{eq3.11}
\varphi(t)\cdot\Bigg(\frac{1}{1 - \frac{t}{v_{(1)}^+}}\Bigg)^{\Big(\sum\limits_{i\in I_1^+}\alpha_i - \sum\limits_{r\in \tilde{I}_1^+}\alpha_r\Big)}\cdots\Bigg(\frac{1}{1 - \frac{t}{v_{(d_+)}^+}}\Bigg)^{\Big(\sum\limits_{i\in I_{(d_+)}^+}\alpha_i - \sum\limits_{r\in \tilde{I}_{(d_+)}^+}\alpha_r\Big)}\\\cdot\Bigg(\frac{1}{1 + \frac{t}{v_{(1)}^-}}\Bigg)^{\Big(\sum\limits_{j\in I_1^-}\alpha_j - \sum\limits_{s\in \tilde{I}_1^-}\alpha_s\Big)}\cdots\Bigg(\frac{1}{1 + \frac{t}{v_{(d_-)}^-}}\Bigg)^{\Big(\sum\limits_{j\in I_{(d_-)}^-}\alpha_j - \sum\limits_{s\in \tilde{I}_{(d_-)}^-}\alpha_s\Big)} = 1
\end{multline}
for $t\neq 0$ and $t\in\big(-v_{(1)}^-, v_{(1)}^+\big)$. 
Observe now that if $\sum\limits_{i\in I_1^+}\alpha_i \neq \sum\limits_{r\in \tilde{I}_1^+}\alpha_r$, then
\[\lim\limits_{t\nearrow v_{(1)}^+}\Bigg(\frac{1}{1 - \frac{t}{v_{(1)}^+}}\Bigg)^{\Big(\sum\limits_{i\in I_1^+}\alpha_i - \sum\limits_{r\in \tilde{I}_1^+}\alpha_r\Big)}\in \{0, +\infty\}\]
which is a contradiction to \eqref{eq3.11}. Therefore, we have
\[\sum\limits_{i\in I_1^+}\alpha_i = \sum\limits_{r\in \tilde{I}_1^+}\alpha_r\]
which by assumption implies that $I_1^+ = \tilde{I}_1^+$. This implies that for all $i\in I_1^+$, we have
\[\frac{\lambda_i}{\beta_i} = v_{(1)}^+ = \frac{\lambda_i}{\tilde{\beta}_i}\]
and therefore
\[\beta_i = \tilde{\beta}_i.\]
The argument can be extended to $I_j^+$ for all $j\in\{1,\ldots,d_+\}$ and to $I_j^-$ for all $j\in\{1,\ldots,d_-\}$ and we therefore get that $\beta_i = \tilde{\beta}_i$ for all $i\in\{1,\ldots,d-1\}$. Equation \eqref{eq3.11} then gives us
\[\varphi(t) = 1,\quad\text{for all }t\in\R\]
which is true if and only if $\mu (\beta_d - \tilde{\beta}_d)  t +  \sigma^2/2 (\beta^2_d  - \tilde{\beta}^2_d) t^2 =0 $  for all $t \in \mathbb R$, and hence $\tilde{\beta}_d  = \beta_d$ since $\mu \ne 0$. \hfill \(\square\)

\subsubsection*{Proof of Lemma \ref{lem:centering_scaling}}
We only prove the inclusion $\subseteq$; the second equality is obvious.

Let $(a, b)$ be an element of $\{(a, b) \in \R^2: X_1 a + X_2 b  \overset{d}{=}  X_1 \alpha + X_2 \beta \}$. It follows that
\begin{align*}
    &X_1 a + X_2 b  - \alpha \mu_1 - \beta \mu_2 \overset{d}{=}  X_1 \alpha + X_2 \beta - \alpha \mu_1 - \beta \mu_2 \\
\Longrightarrow \quad &  \wt{X}_1 a + \wt{X}_2 b + (a  - \alpha) \mu_1 + (b  - \beta) \mu_2 \overset{d}{=}  \wt{X}_1 \alpha + \wt{X}_2 \beta \\
\Longrightarrow \quad& \wt{X}_1 a + \wt{X}_2 b \overset{d}{=}  \wt{X}_1 \alpha + \wt{X}_2 \beta, 
\end{align*}
which yields the assertion. Note that the second implication follows from the fact that
$(a  - \alpha) \mu_1 + (b  - \beta) \mu_2 = 0$ since otherwise the left and right hand side would have different means. \qed

\subsubsection*{Proof of Theorem \ref{prop: nonIID4thmom}}
Note that for any $a,b$ such that $Y_0 \overset{d}{=} a X_1 + b X_2$, it must hold that
\begin{equation*}
\E[(a X_1 + b X_2)^k] = \E[Y_0^k], \quad k=1,\ldots,4. 
\end{equation*}
For $k = 2$, this condition yields
\begin{equation}\label{eq:cond2}
a^2 + b^2  = \alpha^2 + \beta^2 \invcoloneq r^2,
\end{equation}
where $r$ denotes the radius of the sphere containing $\mc{B}_0$. For $k = 4$, the condition yields
\begin{align}
& a^4 \E[X_1^4] + 6 a^2 b^2 \E[X_1^2 X_2^2] + b^4 \E[X_2^4] = \E[Y_0^4] = c r^4 \notag \\
& \; \Longleftrightarrow \ a^4 m_1 + 6 a^2 b^2 + b^4 m_2 = c r^4. \label{eq:cond4} 
\end{align}
Squaring \eqref{eq:cond2} yields
\begin{equation*}
a^4 + b^4 + 2 a^2 b^2 = r^4
\end{equation*}  
and hence $6 a^2 b^2 = 3 r^4 - 3 a^4 - 3b^4$. Back-substituting this into \eqref{eq:cond4}, we obtain
\begin{equation}\label{eq:cond4_back}
a^4 (m_1 - 3) + b^4 (m_2 - 3) = (c -3)r^4.
\end{equation}
Assuming for now that $c \neq 3$ (which implies that not both of $m_1, m_2$ are equal to $3$), we have 
\begin{equation}\label{eq:quartic}
a_u^4 \frac{m_1 - 3}{c - 3} + b_u^4 \frac{m_2 - 3}{c-3} = 1,
\end{equation}  
where $a_u = a/r$ and $b_u = b/r$. Letting $0 \leq x = a_u^2 \leq 1$ so that $0 \leq 1 -x = b_u^2 \leq 1$ the previous equation can be re-written as 
\begin{equation}\label{eq:quadratic}
x^2 \frac{m_1 - 3}{c - 3} + (1 - x)^2 \frac{m_2 - 3}{c-3} = 1. 
\end{equation}
Note that the solutions of the system of equations $\{ x = a_u^2, \ b_u^2 = 1 - x  \}$ are 
\begin{equation*}
  \left \{ (a_u, b_u) \in \{ (\sqrt{x}, \sqrt{1 - x}), (-\sqrt{x}, \sqrt{1 - x}), (\sqrt{x}, -\sqrt{1 - x}), (-\sqrt{x}, -\sqrt{1 - x}) \right  \}.
\end{equation*}
Since the quadratic equation \eqref{eq:quadratic} admits at most two solutions, we conclude that the set of all solutions is at most $ 4 \times 2 = 8$. 
As  defined in the statement of the theorem, let 
\begin{equation*}
w_1 = \frac{m_1 - 3}{c - 3}, \ \text{and} \   w_2 = \frac{m_2 - 3}{c-3}.  
\end{equation*}
Then,  \eqref{eq:quadratic} can be  re-expressed as $(w_1 + w_2) x^2 - 2xw_2 + w_2 - 1 = 0$. Note that if $w_1 + w_2 = 0$, then $w_2 \neq 0$  because otherwise $m_1 = m_2 = 3$, and hence the quadratic equation reduces to a linear one with the single solution $(w_2 - 1)/(2w_2)$\footnote{Note that $(w_2 - 1)/(2w_2) =  (1/2) (1 - 1/w_2)  \in [0,1]$ since $1/\vert w_2 \vert \le 1$. In fact, it follows from (\ref{eq:quartic}) and $w_1 + w_2 = 0$ that $(\beta^4_u  - \alpha_u^4)  w_2 = 1$, where $\alpha_u =  \alpha/r, \beta_u = \beta/r$. Hence, $\vert w_2 \vert \ge 1$ since $\vert \beta^4_u  - \alpha_u^4 \vert \le 1$}.  We conclude that if $w_1 + w_2 = 0$,  the possible solutions belong to a set of cardinality equal to 4. In other words, $(\alpha, \beta)$ is identifiable up to a sign a flip,  and the assertion of the theorem follows.

Now, we consider the case where $w_1 + w_2 \neq 0$. Then the roots of the quadratic equation \eqref{eq:quadratic} are given by
\begin{equation}\label{eq:roots}
\frac{w_2}{w_1 + w_2} \pm \frac{\sqrt{w_1 + w_2 - w_1 w_2}}{w_1 + w_2}.
\end{equation}
Note that $w_1 + w_2 \ge w_1 w_2$; cf.~Lemma \ref{lem:prop_weights} at the end of this appendix. The same lemma shows that it always hold that $\max(w_1, w_2) > 0$. By symmetry, we may assume  w.l.o.g.  that $w_1 > 0$ and $w_1 \geq w_2$.  Suppose that $w_2 < 0$. We now distinguish the following cases.  
\begin{itemize}
\item If $w_1 + w_2 > 0$, then it the root 
\begin{equation*}
\frac{w_2}{w_1 + w_2} - \frac{\sqrt{w_1 + w_2 - w_1 w_2}}{w_1 + w_2}  < 0
\end{equation*}
and hence is outside the target interval $[0,1]$ (recall that $0 \leq x \leq 1$).
\item If $w_1 + w_2 < 0$, then 
      \begin{eqnarray*}
        \frac{w_2}{w_1 + w_2} + \frac{\sqrt{w_1 + w_2 - w_1 w_2}}{w_1 + w_2}  - 1  =  \frac{w_1}{w_1 + w_2} + \frac{\sqrt{w_1 + w_2 - w_1 w_2}}{w_1 + w_2}  \ge 0
      \end{eqnarray*}
      which means that one of the roots (\ref{eq:roots}) is outside $[0,1]$.  
\end{itemize}
Now, note that in either case we have that $\min\{w_1, w_2\}  = w_2 < 0  < 1$, the claim of the proposition holds true. 

Finally, we consider the case
$w_2 \geq 0$ in which case $w_1 + w_2 > 0$. From \eqref{eq:roots}, observe that the existence of two roots inside $[0,1]$ entails that the following two conditions hold:
\begin{equation*}
\sqrt{w_1 + w_2 - w_1 w_2} \leq w_1, \quad  \text{and}  \   \ \    \sqrt{w_1 + w_2 - w_1 w_2} \leq w_2.
\end{equation*}
Squaring both sides of the conditions yields 
\begin{equation*}
(w_1 + w_2) (1 - w_2) \leq 0, \quad  \text{and}  \   \ \    (w_1 + w_2) (1 - w_1) \leq 0.
\end{equation*}
Since $w_1 + w_2 > 0$,  the conditions are equivalent to 
\begin{equation*}
w_1 \ge 1, \quad \text{and} \  \ \  w_2 \ge 1.
\end{equation*}
that is to $w_2 \ge 1$ since we assumed $w_1 \geq w_2$ w.l.o.g. that $w_2 \geq 1$.   We conclude that in the case $w_2 \ge 0$, one of the roots (\ref{eq:roots}) is outside $[0,1]$ if $w_2 = \min(w_1, w_2) < 1$. This means that the condition $\min(w_1, w_2) < 1$ implies again that $(\alpha, \beta)$ is identifiable up to a sign flip.

We now turn to the case $c = 3$. From \eqref{eq:cond4_back}, we deduce -- in the same manner as above -- the   quadratic equation
\begin{equation}\label{eq:quadratic_2}
x^2 \Delta_1 + (1 - x)^2 \Delta_2 = 0, \quad \Delta_1 \coloneq m_1 - 3, \; \Delta_2 \coloneq m_2-3. 
\end{equation}
If $\Delta_1 + \Delta_2 =  0$, then $(1 - 2 x) \Delta_2 = 0$ and hence either $\Delta_2 = 0$ or $x = 1/2$. The first case is not possible because not both $m_1$ and $m_2$ are equal to $3$. Thus, we must have $x  = 1/2$ and $(\alpha, \beta)$ is identifiable to a sign flip. 

Suppose that $\Delta_1 + \Delta_2 \ne 0$. Then, the roots of the equation (\ref{eq:quadratic_2}) are calculated as
\begin{equation*}
\frac{\Delta_2}{\Delta_1 + \Delta_2} \pm \frac{\sqrt{-\Delta_1 \Delta_2}}{\Delta_1 + \Delta_2}. 
\end{equation*}  
Observe that $c = 3$ implies that $\Delta_1,\Delta_2$ have opposite signs. W.l.o.g., suppose $\Delta_1 > 0$ and $\Delta_2 < 0$. Note that
if $\Delta_1 + \Delta_2 > 0$, then one of the roots will be negative (i.e., outside $[0,1]$). Similarly, if $\Delta_1 + \Delta_2 < 0$, then
one of the two roots will exceed 1. $\qed$

\subsubsection*{Auxiliary Lemma}

\begin{lemma}\label{lem:prop_weights} 
Let $w_1$ and $w_2$ be defined as in Theorem \ref{prop: nonIID4thmom}.  Then $w_1 + w_2 - w_1 w_2 \geq 0$. In particular, $\max\{w_1, w_2\} > 0$.
Moreover, if $w_1 = w_2 = w$, then $1 \leq w \leq 2$. 
\end{lemma}

\begin{bew}
We have 
\begin{equation*}
w_1 \cdot w_2 = \frac{(m_1 - 3) (m_2 - 3)}{(c-3)^2}
\end{equation*}
and 
\begin{align*}
w_1 + w_2 = \frac{(m_1 - 3) + (m_2 - 3)}{c-3} &=  \frac{\{(m_1 - 3) + (m_2 - 3)\} (c - 3)}{(c-3)^2} \\
 &= \frac{\{(m_1 - 3) + (m_2 - 3)\} \{ (m_1 - 3) \alpha_u^4 + (m_2 - 3) \beta_u^4 \}}{(c-3)^2}.  
\end{align*}
We also have 
\begin{align*}
(m_1-3)(m_2-3) &= (m_1-3)(m_2-3)(\alpha_u^2 + \beta_u^2)^2 \\
               &= (m_1 - 3)(m_2-3)(\alpha_u^4 + \beta_u^4) + 2\alpha_u^2\beta_u^2 (m_1-3)(m_2-3) \\
               &\leq (m_1 - 3)(m_2-3)(\alpha_u^4 + \beta_u^4) + (m_1-3)^2 \alpha_u^4 + (m_2-3)^2 \beta_u^4,
\end{align*}
which yields the claim $w_1 + w_2 - w_1 w_2 \geq 0$. Regarding the ``Moreover, $\ldots$'' part, observe that
$w_1 = w_2$ implies that $m_1 = m_2 = m$ and thus $c - 3 = \alpha_u^4 m + \beta_u^4 m - 3(\alpha_u^2 + \beta_u^2)^2$, which
yields $c - 3 = (m-3)(\alpha_u^4 + \beta_u^4)$ and in turn $w = 1/(\alpha_u^4  + \beta_u^4)$, so that $1 \leq w \leq 2$.
\end{bew}

\end{document}